\documentclass{article}
\usepackage{amsmath}
\usepackage{amssymb}

\usepackage{graphicx}

\def\av{\textrm{Av}}

\def\degree{\textrm{degree}}

\def\cat{\mbox{\rm cat}}
\def\Nat{\mathbb{N}}

\newcommand{\W}[2]{W_{#1 \! \makebox[1pt]{} #2}}

\renewcommand{\int}{\cap}

\newcommand{\intersection}{\cap}

\newtheorem{lemma}{Lemma}
\newtheorem{theorem}[lemma]{Theorem}

\newtheorem{proposition}[lemma]{Proposition}
\newtheorem{corollary}[lemma]{Corollary}

\newcommand{\qed}{\hfill \rule{1ex}{1ex}} 
\newenvironment{proof}{{\bf Proof}: }{\qed}

\title{Permutation Classes of Polynomial Growth}
\author{
M.~H. Albert\thanks{{\tt malbert@cs.otago.ac.nz}},
M.~D. Atkinson\thanks{{\tt mike@cs.otago.ac.nz}}\\
Department of Computer Science\\
University of Otago, Dunedin, New Zealand\\
and\\
Robert Brignall\thanks{{\tt robertb@mcs.st-and.ac.uk}}\\
Department of Mathematics and Statistics\\
St Andrews, UK.}

\begin{document}
\maketitle
\begin{abstract}
A pattern class is a set of permutations closed under the formation
of subpermutations. Such classes can be characterised as those
permutations not involving a particular set of forbidden
permutations. A simple collection of necessary and sufficient
conditions on sets of forbidden permutations which ensure that the
associated pattern class is of polynomial growth is determined. A
catalogue of all such sets of forbidden permutations having three or
fewer elements is provided together with bounds on the degrees of
the associated enumerating polynomials.
\end{abstract}
\section{Introduction}
A permutation $\pi$ is said to be a subpermutation of a permutation $\sigma$, $\pi\preceq\sigma$, if $\sigma$ has a subsequence isomorphic to $\pi$ (that is, its terms are ordered relatively the same as the terms of $\pi$).  For example $312$ is a subpermutation of $25134$ because of the subsequence $513$ (or $514$ or $534$). On the other hand $321$ is not a subpermutation of $25134$ because there is no three element subsequence of $25134$ in which  the three elements occur in decreasing order. Consequently
$25134$ is said to {\em involve}\ $312$ but to {\em avoid}\ $321$. If $\Pi$ is a set of permutations then $\av(\Pi)$ denotes the set of all permutations which avoid every permutation in $\Pi$.  Such sets of permutations are called \emph{pattern classes}\ and have given rise to many enumerative
results.  Typically, given $\Pi$, one is interested in determining the number $c_n(\Pi)$ of permutations of each length $n$ in the pattern class $\av(\Pi)$. For obvious reasons we shall assume throughout that $\Pi$ is non empty. When explicitly listing the elements of some set $\Pi$ as an argument we will generally omit braces, thus writing $c_n(123, 312)$ rather than $c_n(\{123, 312\})$.

The sequences $c_n(\Pi)$  can be studied from several points of view.  We might wish to discover an exact formula for $c_n(\Pi)$, to find bounds on its growth as a function of $n$, or to determine the ordinary generating function
\[
\sum_{\sigma\in\av(\Pi)}x^{|\sigma|}.
\]
Recently Marcus and Tardos \cite{Marcus} resolved affirmatively the long-standing open question of whether $c_n(\Pi)$ was always exponentially bounded.  In part because of this result, attention has turned to enumerative questions of finer detail, and this paper addresses one such question.

We shall be concerned with pattern classes of polynomial growth; in other words, classes $\av(\Pi)$ for which there exists a bound of the form
\[
c_n(\Pi)\leq An^d
\]
for some constants $A, d$. Examples of classes of polynomial growth have appeared many times in the literature.  For example, in an early paper \cite{SS}, on pattern class enumeration, Simion and Schmidt proved that
$c_n(132, 321)=n(n-1)/2+1$.  Some more difficult enumerations were carried out by West \cite{West} in his work on classes of the form $\av(\alpha,\beta)$ where $\alpha$ is a permutation of length three, and $\beta$ one of length four; he showed that 4 of the 18 essentially different such classes have polynomial enumerations.

More recently, Kaiser and Klazar \cite{KK} proved that, in a polynomial growth class, $c_n(\Pi)$, as a function of $n$, was actually equal to some polynomial for all sufficiently large $n$ and that this polynomial had a particular form. Kaiser and Klazar also proved that classes $\av(\Pi)$ whose growth was not polynomial have $c_n(\Pi)\geq \tau^n$ where $\tau$ is the golden ratio.

Huczynska and Vatter \cite{HV} gave a simplification of the results of \cite{KK}, characterising polynomial growth classes in terms of ``grid classes'' of matchings, and establishing the dichotomy between classes of polynomial growth and those whose growth exceeds the growth of the Fibonacci numbers.

Necessary and sufficient conditions on the basis $\Pi$ for
$\av(\Pi)$ to have polynomial growth are implicit in \cite{KK}, and
were made explicit by Huczynska and Vatter. These conditions  are so
simple that it is virtually trivial to test whether $\av(\Pi)$ has
polynomial growth. By themselves the conditions tell us little about
an actual polynomial that gives $c_n(\Pi)$ (for sufficiently large
$n$) and so, after exhibiting a somewhat different derivation of the
conditions (Theorem \ref{N&S}), we go on to give more precise
results when $|\Pi|\leq 3$ (Section \ref{small}).

If $\Pi = \{ \alpha \}$ there is nothing to say beyond what is obvious; $c_n(\alpha)$ has polynomial growth only if $|\alpha| \leq 2$. In these cases:
\begin{eqnarray*}
c_n(1) &=& 0 \quad \mbox{for all $n \geq 1$} \\
c_n(12) = c_n(21) &=& 1\quad \mbox{for all $n \geq 1$}
\end{eqnarray*}

It is more complex to characterise the two and three element bases $\Pi$ that lead to pattern classes of polynomial growth. In the three element case the classes $\av(\Pi)$ of polynomial growth are sufficiently numerous that we have only used Theorem \ref{N&S} to list the various sets $\Pi$ (see Theorem \ref{threerestrictions}); it would not be difficult in most cases to give the complete enumerations.  However, in the two element case, we obtain a characterisation (Theorem \ref{tworestrictions}) of polynomial growth classes which are more demanding to analyse.  In  Section \ref{enumeration}  we give some bounds on the degrees of the polynomials that arise in this case.

In order to simplify the exposition it will be useful to introduce a few further pieces of definition and notation. Two sequences $a_1, a_2, \ldots, a_n$ and $b_1, b_2, \ldots, b_n$ of distinct elements from (possibly different) totally ordered sets are {\em order isomorphic}\ (or simply {\em equivalent}) if, for all $1 \leq i, j \leq n$, $a_i < a_j$ if and only if $b_i < b_j$. Thus, a permutation $\pi$ is involved in a permutation $\sigma$ exactly when, considered as a sequence, it is equivalent to some subsequence of $\sigma$. Further, every finite sequence of distinct elements from a totally ordered set is equivalent to exactly one permutation, called its {\em pattern}. If a pattern class $X = \av(\Pi)$ is of polynomial growth, then we define $\degree(X)$ to be the degree of the polynomial $p$ for which $c_n(\Pi) = p(n)$ for all sufficiently large $n$.

\section{Conditions for polynomial growth}
\begin{figure}[ht]
\begin{center}
\includegraphics[width=3in]{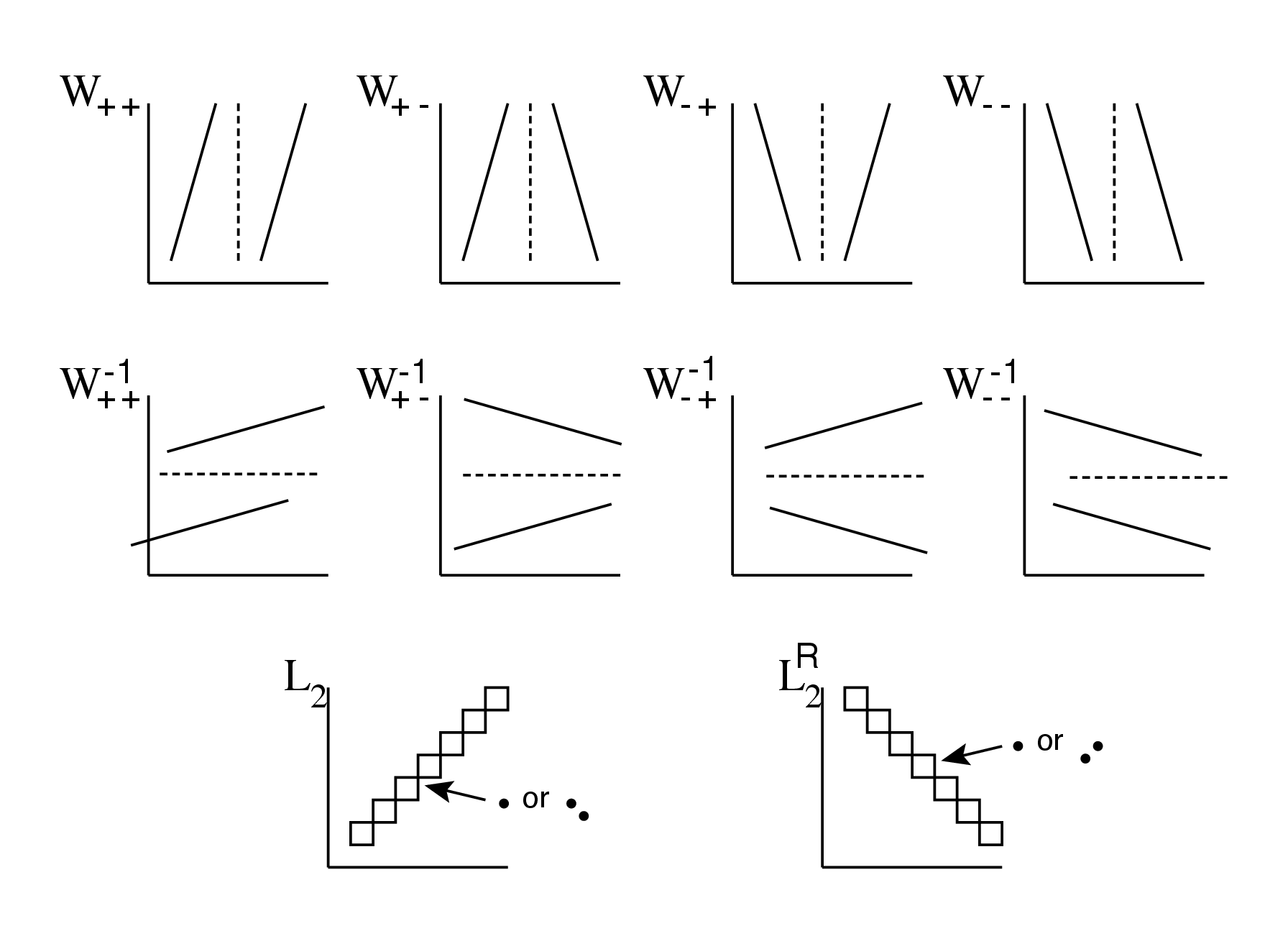}
\caption{10 types of permutation}
\label{10types}
\end{center}
\end{figure}

This section aims to give an alternative proof to \cite{HV} of the
necessary and sufficient condition for which $\av(\Pi)$ has
polynomial growth.  Informally this condition is that among the
permutations of $\Pi$ we must find permutations of all the 10 types
shown in Figure \ref{10types}. Clearly, testing this condition is
very easy and can be done in linear time.

We shall develop some terminology and notation to state this condition more formally, and to justify it.  Let $\epsilon=(e_1,e_2,\ldots,e_r)$ be any sequence whose terms are $+1$ or $-1$.  Then the pattern class
$W(\epsilon)$ consists of all permutations $\pi$ that have a segmentation
\[\pi=\sigma_1\sigma_2\cdots\sigma_r\]
where $\sigma_i$ is increasing if $e_i=+1$ and decreasing if $e_i=-1$ (empty segments are allowed). For example, the permutation $\pi=35178426\in W(+1,+1,-1,+1)$ because of the segmentation $35|178|42|6$. These pattern classes are the `$W$'-classes of \cite{pwo,regular} where they were used to study partial well-order and regularity questions.  We will be particularly interested in the four $W$-classes formed from sequences $\epsilon$ of length $2$ and their inverses, which form the first eight types in Figure \ref{10types}.  We will use a somewhat more compact notation for these classes:
\[
\W{+}{-} = W(+1, -1) \quad \W{-}{-}^{-1} = W(-1,-1)^{-1} \quad \mbox{etc.}
\]
The last two types in Figure \ref{10types} are not related to $W$-classes and we call them $L_2$ and $L_2^R$ respectively since the first consists of permutations with increasing layers which are either singletons or decreasing doubletons, and the second is the reverse of this class.  We can now state the conditions formally:

\begin{theorem} \label{N&S}A pattern class $X=\av(\Pi)$ has polynomial growth if and only if every class in the list
\[
\begin{array}{c}
\W{+}{+} , \, \W{+}{-} , \, \W{-}{+} , \, \W{-}{-}  \\
\W{+}{+}^{-1} , \, \W{+}{-}^{-1} , \, \W{-}{+}^{-1} , \, \W{-}{-}^{-1}  \\
L_2, \, L_2^R
\end{array}
\]
 has non-empty intersection with  $\Pi$.
\end{theorem}

Note that in order to avoid each of the ten types, we may assume
that $\av(\Pi)$ has at most $10$ basis elements. Our proof of
Theorem~\ref{N&S} needs two auxiliary results.

\begin{lemma}Suppose that $X$ is a pattern class that is not a subset of any
$W$-class.  Then there exists a one-to-one map
 $\gamma:\Nat\longrightarrow\mathbb{R}$ with
\[
\gamma(1)<\gamma(2)>\gamma(3)<\gamma(4)>\gamma(5)<\ldots
\]
such that $X$ contains every subpermutation of $\gamma$.
\end{lemma}
\begin{proof}
A permutation $\pi$ is called an {\em up-down}\ permutation if $\pi_1 < \pi_2 > \pi_3 < \pi_4 \cdots $. The set of all up-down permutations can be considered as a tree, $T$. The root of this tree is the single permutation of length $1$ and the $k$th level is the set of up-down permutations of length $k$.  The parent of a permutation on the $k$th level is the
up-down permutation of length $k-1$ defined by removing its last symbol and replacing the remaining elements by their pattern. For example, the parent of $2 6 1 5 3 4$ is $2 5 1 4 3$. In particular, note that $T$ is finitely branching, with each vertex at level $k$ having at most $k$ children.

Since $X$ is not contained in any $W$ class, it must contain a basis
element of each $W$ class. Among these there are up-down
permutations of every length. These form an infinite subtree of $T$
as $X$ is closed under taking subpermutations. By K\"{o}nig's
lemma\footnote{\textbf{K\"onig's Lemma.} A finitely branching tree
is infinite if and only if it has an infinite path.}, this tree has
an infinite path.  Thus $X$ contains a sequence of up-down
permutations $\delta_1, \, \delta_2, \ldots$ such that, for each
$i$, the length of $\delta_i$ is $i$ and $\delta_{i+1}$ is a child
of $\delta_i$.

By a general construction in \cite{Natural} (Theorem 1.2), the set
of subpermutations of this sequence of permutations can be described
in terms of a map $\gamma$ where $\gamma$ has the form claimed in
the statement of the lemma. In this case it is easy to carry out the
construction directly: Inductively define the map $\gamma$ (starting
from an empty map) by defining $\gamma(i+1)$ to be any real number
such that the sequence $\gamma(1) \gamma(2) \cdots \gamma(i+1)$ is
equivalent to $\delta_{i+1}$.
\end{proof}

\begin{proposition}\label{minimal-non-W}
Suppose that the pattern class $X$ is not contained in any $W$-class.  Then $X$ contains one of the following pattern classes:
$\W{+}{+}^{-1}$, $\W{+}{-}^{-1}$, $\W{-}{+}^{-1}$, $\W{-}{-}^{-1}$,
$L_2$, or
$L_2^R$.
\end{proposition}
\begin{proof}
Consider the map $\gamma$ guaranteed by the previous lemma. For
$n\geq 1$, the sequence of real numbers $\gamma(2n)$ contains an
infinite monotone subsequence, labelled
$\gamma(2n_1),\gamma(2n_2)\ldots$.

Suppose first that this sequence is increasing, and put
$a_{2i}=\gamma(2n_i)$, $a_{2i+1}=\gamma(2n_i+1)$, noting that
$a_{2i}>a_{2i+1}$ for all $i$. Now consider the infinite complete
graph with vertices $r$ for $r$ a positive integer. We colour the
edges $(r,s)$ of this graph with the pattern of the sequence $a_{2r}
a_{2r+1} a_{2s} a_{2s+1}$. Given the constraints $a_{2r} >
a_{2r+1}$, $a_{2r} < a_{2s}$ and $a_{2s} > a_{2s+1}$ there are only
three possibilities for this pattern: $2143$, $3142$ and $3241$. By
Ramsey's theorem the graph contains an infinite monochromatic
induced subgraph. Since all permutations equivalent to subsequences
of the values of $\gamma$ belong to $X$ there is no loss of
generality in assuming that the entire graph is monochromatic.
However, this implies that $X$ contains $L_2$, $\W{+}{+}^{-1}$ or
$\W{-}{+}^{-1}$ depending on whether the colour occurring is $2143$,
$3142$ or $3241$.

For the other case, when $\gamma(2n)$ contains an infinite
decreasing sequence, put $a_{2i-1}=\gamma(2n_i-1)$,
$a_{2i}=\gamma(2n_i)$ so that $a_{2i-1}<a_{2i}$, and consider
instead the graph whose edges $(r,s)$ are coloured with the pattern
of the sequence $a_{2r-1} a_{2r} a_{2s-1} a_{2s}$. Then an exactly
parallel argument to that of the previous paragraph establishes that
$X$ must contain one of $L_2^R$, $\W{-}{-}^{-1}$ or $\W{+}{-}^{-1}$.
Alternatively we could note that in this case the class consisting
of the reversals of all the permutations in $X$ must be of the type
already analysed.
\end{proof}

We can now complete the proof of Theorem \ref{N&S}.

\begin{proof}
One implication is clear: if $X$ has polynomial growth then it
cannot contain any of the 10 classes specified since these all have
exponential growth, and hence $\Pi$ must contain a permutation from
each of them.

For the converse we shall use Proposition \ref{minimal-non-W}.  So
now suppose that $\Pi$ contains a permutation from each of the 10
given pattern classes.  Then $X$ does not contain any of the 10
classes.  Then, by Proposition \ref{minimal-non-W}, it must be
contained in some $W$-class $W_{1 }$ (with $a$ segments, say).  But,
applying Proposition \ref{minimal-non-W} to $X^{-1}$, it must also
be the case that $X$ is contained in the inverse of some $W$-class
$W_{2}$ (with $b$ segments say).  Consider a permutation $\pi$
belonging to the intersection of $W_1$ and $W_2^{-1}$. Its
representation as an element of $W_1$ divides the positions of $\pi$
into $a$ blocks, in each of which the values form a monotone
segment. Similarly, its representation as an element of $W_2^{-1}$
divides the values of $\pi$ into $b$ blocks of consecutive elements.
These value blocks might cross-cut each of the $a$ position blocks,
and likewise the position blocks might cross-cut the value blocks.
However, $\pi$ will have a blocked structure with (at most)
$ab\times ab$ blocks, where the $ab$ non-empty blocks are monotone
and form a permutation pattern. Figure \ref{intersect} provides a
simple illustration of this.

Since the number of such permutations of length $n$ is bounded above
by the number of non-negative solutions of $x_1 + x_2 + \cdots +
x_{ab} = n$, which is a polynomial in $n$ of degree $ab-1$, $X$ has
polynomial growth.
\end{proof}

This proof shows that every class of polynomial growth is a subclass
of a polynomial growth class defined by a permutation $\pi$ (of
degree $m$ say) and a sequence of $m$ signs $\pm 1$.  To avoid a
reduction to a smaller case we will assume that when this
permutation has consecutive terms $i,i+1$ then not both signs are
$+1$, and when it has consecutive terms $i+1,i$ then not both signs
are $-1$.  The permutations in the class are obtained from $\pi$ by
replacing any term associated with $+1$ by an increasing consecutive
segment (possibly empty), and the terms associated with $-1$ by a
decreasing consecutive segment.  Therefore any permutation in the
class can be specified (though not generally uniquely) by the vector
of lengths of these segments.  A subclass  then corresponds to an
ideal in the partially ordered set of such vectors ordered by
dominance.

\section{Two or three restrictions}\label{small}

In this section we consider the implications of Theorem \ref{N&S} for $\Pi$ when $|\Pi|=2$ or $3$. To eliminate trivialities we will assume throughout this section that each permutation in $\Pi$ has length at least three.

\begin{theorem}\label{tworestrictions}
The class $X=\av(\alpha,\beta)$ has polynomial growth if and only if
(up to symmetry and exchange of $a$ with $b$) we have one of the
following:
\begin{enumerate}
\item $\alpha$ is increasing and $\beta$ is decreasing,
\item $\alpha$ is increasing and $\beta$ is almost decreasing in the sense that $\beta\in L_2^R$ with exactly one layer of size $2$.
\end{enumerate}
\end{theorem}
\begin{proof}

One implication is obvious; it is clear that each of the stated
classes has polynomial growth by Theorem~\ref{N&S}.

On the other hand, since $\{\alpha,\beta\}$ contains at least one
element from $L_2$ we may, by exchanging $\alpha$ with $\beta$ if
necessary, assume that $\alpha\in L_2$.  But also one of
$\alpha,\beta$ is in $L_2^R$. However $\alpha\in L_2^R$ leads to a
contradiction since $L_2\cap L_2^R=\{1,12,21\}$ and $|\alpha|>2$, so
from now on we assume that $\beta\in L_2^R$.

Next we shall prove that one of $\alpha,\beta$ is monotone.  Assume for a contradiction that neither is monotone.  Since $\W{+}{-} \cap \W{-}{+}$ consists only of monotone permutations one of $\alpha,\beta$ belongs to
$\W{+}{-}$ and the other belongs to $\W{-}{+}$.  By the same reasoning one of $\alpha,\beta$ belongs to
$\W{+}{-}^{-1}$ and the other belongs to $\W{-}{+}^{-1}$.  Up to symmetry there are only two possible cases:
\[
\begin{array}{l}
\alpha\in \W{+}{-} \cap \W{+}{-}^{-1}, \: \beta\in \W{-}{+}\cap \W{-}{+}^{-1} \quad \mbox{or} \\
\alpha\in \W{+}{-} \cap \W{-}{+}^{-1}, \: \beta\in \W{-}{+}\cap \W{+}{-}^{-1}.
\end{array}
\]
In the first case the only non-monotone permutations of $\W{-}{+} \cap \W{-}{+}^{-1}$ are of the form $k(k-1)\cdots 1(k+1)(k+2)\cdots n$, but none of them are in $L_2^R$ and so $\beta$ must be monotone, a contradiction.  For the other case a similar contradiction can also be obtained.

We shall assume (using symmetry) that $\alpha$ is increasing.  As $|\alpha|>2$, we have $\alpha\not\in \W{-}{-} \cup \W{-}{-}^{-1}$ and therefore we have $\beta\in \W{-}{-} \cap \W{-}{-}^{-1}$.  This class is easily seen to consist of permutations whose shape is as shown in Figure \ref{intersect}; but $\beta$ also lies in $L_2^R$ and so has the form given in the theorem.
\begin{figure}[ht]
\begin{center}
\includegraphics[width=2in]{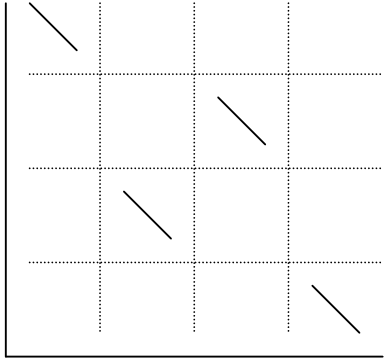}
\end{center}
\caption[]{The shape of permutations in \protect{\mbox{$\W{-}{-} \cap \W{-}{-}^{-1}$}}.}
\label{intersect}
\end{figure}
\end{proof}

\begin{theorem}\label{threerestrictions}
The class $\av(\alpha,\beta,\gamma)$ has polynomial growth if and
only if (up to symmetry and exchange of $\alpha$, $\beta$ and
$\gamma$) we have one of the following:
\begin{enumerate}
\item $\alpha = 213$, and
\begin{enumerate}
\item  $3412\preceq\beta\in L_2^R$ and $\gamma = 12 \cdots k \, n (n-1) \cdots (k+1)$ for some $k$, or
\item $\beta = m(m-1) \cdots (j+2)\, j(j+1) \,(j-1)(j-2) \cdots 1$ and $\gamma = 12 \cdots k \, n(n-1) \cdots (k+1)$ for some $j, k$, or
\item $\beta = m(m-1) \cdots 312$ and $\gamma \in \W{+}{-}$.
\end{enumerate}
\item $\alpha$ is increasing, and
\begin{enumerate}
\item $3412\preceq\beta\in L_2^R$ and $\gamma \in \W{-}{-} \cap \W{-}{-}^{-1}$, or
\item $\av(\alpha,\beta)$ has polynomial growth.
\end{enumerate}
\item $\alpha = 21345 \cdots  s$, and
\begin{enumerate}
\item $3412\preceq\beta\in L_2^R$ and $\gamma = 1n(n-1) \cdots 2$, or
\item $\beta = m(m-1) \cdots (j+2)\, j(j+1) \, (j-1)(j-2) \cdots 1$ and $\gamma = 12 \cdots k \, n (n-1) \cdots (k+1)$ for some $j,k$, or
\item $\beta = m(m-1) \cdots 312$ and $\gamma\in \W{+}{-}$.
\end{enumerate}
\end{enumerate}
\end{theorem}

\begin{proof}
It is clear that each of the stated classes has polynomial growth by
Theorem~\ref{N&S}.

Conversely, since $\av(\alpha,\beta,\gamma)$ has polynomial growth,
$\{\alpha,\beta,\gamma\}$ has non-empty intersection with each of
the 10 classes of Theorem \ref{N&S}. Since $L_2\intersection L_2^R =
\{1,12,21\}$ none of $\alpha,\beta$ or $\gamma$ lies in both. Thus
we will assume that $\alpha\in L_2$, and $\beta\in L_2^R$.

{\bf Case 1.} First suppose $|\alpha|\leq 3$, and that none of $\alpha,\beta$ or $\gamma $ are monotone (which will be covered in case 2). Up to symmetry, therefore, we may assume $\alpha = 213$.
Since $\beta\in L_2^R$, $\beta$ must consist of decreasing layers, each layer being either a singleton or an increasing doubleton.

Suppose first that $\beta$ has at least two layers of size 2 (so $3412\preceq\beta$). Then neither $\alpha$ nor $\beta$ can lie in $\W{+}{-}$ or $\W{+}{-}^{-1}$, and so $\gamma$ must lie in both, giving $\gamma = 12 \cdots k \, n (n-1) \cdots (k+1)$. This is condition (a).

Now suppose $\beta$ has just one layer of size 2. This layer may lie at the left or right hand ends of $\beta$, or it may lie in the middle. If it lies in the middle, i.e. $\beta = m(m-1) \cdots (j+2)\, j(j+1) \, (j-1)(j-2) \cdots 1$ for some $j$, then again neither $\alpha$ nor $\beta$ lies in $\W{+}{-}$ or $\W{+}{-}^{-1}$ hence $\gamma = 12 \cdots k \, n (n-1) \cdots (k+1)$. This is condition (b).

Finally, suppose $\beta$ has a single layer of size 2 at an end. Without loss of generality we may assume this layer is at the right hand end, and that $\beta = m(m-1) \cdots 312$, since the other case may be obtained by the inversion symmetry (which preserves $\alpha=213$). In this case, at least one of $\alpha$ and $\beta$ lies in every class except $\W{+}{-}$, and so $\gamma\in \W{+}{-}$. This is condition (c).

{\bf Case 2.} Now suppose $\alpha$ is increasing (of length at least 3). Here, $\alpha$ lies in all classes except $\W{-}{-}$, $\W{-}{-}^{-1}$ and $L_2^R$.
We assumed that $\beta$ lay in $L_2^R$. Thus $\beta$ consists of decreasing  layers each of which is either a singleton or an increasing doubleton.

If $\beta$ has at least two layers of size 2 (so $3412\preceq\beta$), then $\beta$ does not lie in $\W{-}{-}$ or in  $\W{-}{-}^{-1}$, and hence $\gamma\in \W{-}{-} \cap \W{-}{-}^{-1}$ whose permutations have the form given in Figure \ref{intersect}.

If $\beta$ has just one layer of size 2, then one of $\alpha$ and $\beta$ lies in each of the 10 classes, so $\av(\alpha,\beta)$ has polynomial growth, and there is no restriction on $\gamma$.

{\bf Case 3.} We may now assume that all of $\alpha,\beta,\gamma$ have length at least 4, and that none are monotone. By our assumption, $\alpha\in L_2$ and $\beta\in L_2^R$, so $\alpha$ and $\beta$ must each have at least one layer of size 2.

Suppose first that both $\alpha$ and $\beta$ have at least two layers of size 2. Then neither $\alpha$ nor $\beta$ can lie in any of $\W{+}{-}$, $\W{-}{+}$,  $\W{+}{-}^{-1}$ or $\W{-}{+}^{-1}$, so $\gamma$ must lie in all of these. However, this would  imply that $\gamma$ was monotone, a contradiction.
Thus we may assume, by interchange of $\alpha$ with $\beta$ and the inversion symmetry if necessary, that $\alpha$ has just one layer of size 2.

Suppose $\beta$ has at least two layers of size 2. If $\alpha$ does not lie in at least one of $\W{+}{-}$ and $\W{-}{+}$, then $\gamma$ will again be monotone, since $\beta$  lies in neither class. Thus $\alpha$ lies in one of $\W{+}{-}$ and $\W{-}{+}$, and we may assume without loss of generality (by the reverse complement symmetry) that $\alpha\in \W{-}{+}$, and so $\alpha$ has the form $21345 \cdots s$.

Neither $\alpha$ nor $\beta$ can lie in $\W{+}{-}$, $\W{-}{-}$ or $\W{+}{-}^{-1}$, and so $\gamma\in \W{+}{-} \intersection \W{-}{-} \intersection \W{+}{-}^{-1}$. The non-monotone permutations in this class are all of the form $1 n (n-1) \cdots 2$, giving condition (a).

So now suppose that $\beta$ has just one layer of size 2. Suppose further that neither $\alpha$ nor $\beta$ had their single layers  at an end, i.e. $1324\preceq\alpha$ and $4231\preceq\beta$. Then neither $\alpha$ nor $\beta$  lie in $\W{+}{-}$ or $\W{-}{+}$, and as we have seen before, this would mean that $\gamma$ was monotone, a contradiction.
Thus we may assume that one of $\alpha$ and $\beta$ has its single layer at an end.  Again by symmetry and interchange of $\alpha$ with $\beta$, we may suppose that $\alpha$ has its size 2 layer at its left end.

Now suppose that $\beta$ has its size 2 layer not at an end, so $4231\preceq\beta$. Then neither $\alpha$ nor $\beta$ lie in $\W{+}{-}$ or $\W{+}{-}^{-1}$, so $\gamma$ must lie in $\W{+}{-}\intersection \W{+}{-}^{-1}$, and therefore, as $\gamma$ is not monotone, $\gamma = 12 \cdots k \, n (n-1) \cdots (k+1)$ for some $k$. This is condition (b).

Finally, suppose that $\beta$ does have its size 2 layer at an end. By inverse symmetry (which is $\alpha$-preserving), we may choose this to be the right end, so $\beta = m (m-1) \cdots 312$. Then neither $\alpha$ nor $\beta$ lies in $\W{+}{-}$, so we have $\gamma\in \W{+}{-}$. This is condition (c).
\end{proof}

\section{Enumeration when there are two restrictions}\label{enumeration}

Throughout this section we shall only consider classes $\av(\alpha, \beta)$ defined by two restrictions of the form given in Theorem \ref{tworestrictions}. In the first case, where $\alpha$ is increasing and $\beta$ is decreasing, $\av(\alpha, \beta)$ is finite by the Erd\H{o}s-Szekeres Theorem. So we will consider only the second case, namely that for some positive integer $r$ and non-negative integers $p$ and $q$:
\begin{enumerate}
\item $\alpha = \alpha_r=12\cdots r$, and
\item $\beta = \beta_{pq} = \lambda \,(q+1) \, (q+2) \,  \mu$ where $|\lambda|=p,|\mu|=q$, $\lambda$ is decreasing with consecutive terms all of which are greater than $q+2$, and $\mu$  is decreasing with consecutive terms, all of which are less than $q+1$.  Define $s=|\beta|=p+q+2$.
\end{enumerate}

We shall give upper and lower bounds on $\degree(\av(\alpha_r, \beta_{pq}))$ for arbitrary $r,p,q$, and some tighter bounds in small special cases.  Our techniques depend on a study of permutations that have no segment of the form $i+1,i$; we call such permutations \emph{irreducible}\/ (this being a slight variation of the terminology of \cite{wreath}).

\subsection{Degree bounds}

\begin{lemma}\label{maxdescent}
An irreducible in the class $\av(\alpha_r,\beta_{pq})$, with $p,q$ both non-zero, has decreasing subsequences of length at most $(r-1)(s-2) - 1$. When one of $p$ or $q$ is zero, an irreducible has decreasing subsequences of length at most $(r-1)(s-2)$.

The length of an irreducible in  $\av(\alpha_r,\beta_{pq})$ is at most:
\[
\begin{array}{cl}
(r-1)^2(s-2) - (r-1) & \mbox{if $p > 0$ and $q > 0$,} \\
(r-1)^2(s-2) & \mbox{if $p = 0$ or $q = 0$.}
\end{array}
\]
\end{lemma}

\begin{proof}
Suppose first that $p>0$ and $q>0$.  Let $\pi \in \av(\alpha_r, \beta_{pq})$ be irreducible, and let $\gamma=g_1g_2\cdots g_d$ be a maximal decreasing subsequence of $\pi$.  We shall deduce properties of $\pi$ using its graph, shown in Figure \ref{upperbound} as laid out in $9$ regions defined by $p,q$ and the points of $\gamma$.  These properties will provide a bound on the number of starred points of $\gamma$.

\begin{figure}[ht]
\begin{center}
\includegraphics[width=2in]{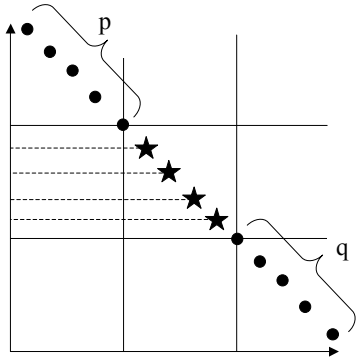}
\caption{An irreducible permutation in $\av(\alpha_r,\beta_{pq})$}
\label{upperbound}
\end{center}
\end{figure}

The middle-left and middle-right regions can be divided into rows
defined by the starred points of $\gamma$. Similarly, the top-middle
and bottom-middle can be divided into columns.  The figure shows the
rows in the middle-left region.

For each non-empty row in the middle-left region choose a representative point of $\pi$.  Among these representatives we cannot have a decreasing sequence $r_1,r_2,\ldots,r_{p+1}$ of length $p+1$.  For if there were such a sequence then we could find a point $u$ among the starred points, below $r_{p}$ and above $r_{p+1}$, and then
\[r_1,r_2,\ldots,r_{p+1},u,g_{d-q+1},\ldots,g_d\]
would be equivalent to $\beta_{pq}$.

An exactly similar argument applies to the representative points of
columns in the top-middle region.  Now consider the representatives
in both the middle-left and top-middle regions.  Among these there
cannot be a decreasing subsequence of length $p+1$.  Nor can there
be an increasing subsequence of length $r-1$ (or together with $g_p$
we would have an increasing subsequence of length $r$, i.e. a
subsequence isomorphic to $\alpha_r$).  By the Erd\H{o}s-Szekeres
theorem the total number of representatives in these two regions is
at most $(r-2)p$.

An analogous argument shows that, in the middle-right and bottom-middle regions, there are at most $(r-2)q$ row or column representatives.  So the total number of non-empty rows or columns in these four regions is at most $(r-2)(p+q)$.

Next we consider the $d-p-q+1$ pairs $(g_p,g_{p+1}), (g_{p+1},g_{p+2}),\ldots,(g_{d-q},g_{d-q+1})$.  None of these pairs can be two consecutive values as $\pi$ is irreducible, and so each of them must have some separating term (in the top-middle or bottom-middle if the separation is by position, and in the middle-left or middle-right if the separation is by value).  But these separating terms clearly all lie in distinct rows or distinct columns, and it follows that
\[
d-p-q+1\leq (r-2)(p+q)
\] which gives
$d\leq (r-1)(s-2)-1$ as required.

Now consider the case where $p=0$. The same argument as above can be used except that here we have $d-p-q$ pairs of terms (rather than $d-p-q+1)$ that require separating representatives; this leads to the required bound. The case $q=0$ follows in a similar way.

The second part of the lemma follows from the Erd\H{o}s-Szekeres theorem since we have a bound on the length of a maximum decreasing subsequence of an irreducible in $\av(\alpha_r, \beta_{pq})$, while an increasing subsequence can have length at most $r-1$.
 \end{proof}

We now discuss how the irreducibles and their properties determine the degree of a polynomial growth class $X$.  Every permutation in $X$ can be contracted to  a unique irreducible permutation by replacing all occurrences of segments $i+1,i$ by $i$ (and relabelling), and doing this repeatedly until the result is irreducible.  For example, $21654873$ reduces to $1342$.  The opposite process of replacing each term in an irreducible permutation of $X$ by decreasing consecutive segments is called {\em expanding}.  Since every permutation in $X$ arises by expanding some irreducible, and since each irreducible of length $m$ expands to $\binom{n-1}{m-1}$ permutations of length $n$ we see that $\degree(X)\leq m-1$ where $m$ is the maximal length of an irreducible permutation in $X$.

However, it does not follow that $\degree(X)= m-1$ since not every expansion of an irreducible lies in $X$.  We therefore introduce the idea of an expansible set: a subset of the terms of an irreducible permutation $\theta$ is said to be \emph{expansible} if the terms can each be replaced simultaneously by an arbitrarily long decreasing consecutive segment to obtain a permutation in the class.  If $\theta$ has an expansible set of size $e$ then certainly $\degree(X)\geq e-1$.

Conversely, if $\degree(X)=e-1$, then $X$ will contain some irreducible permutation which has an expansible subset of size $e$.  To see this note first that $X$ is the finite union of subsets $X_{\pi}$, one for each irreducible $\pi\in X$, where the permutations in $X_{\pi}$ all contract to $\pi$.  Therefore there must be some irreducible $\pi$ for which $X_{\pi}$ has polynomial growth of degree $e-1$.  Suppose that $|\pi|=m$.  Then the permutations of length $n$ in $X_{\pi}$ are determined by a set $Y_n$ of $m$-tuples of positive integers that sum to $n$.  We wish to show that there is some set of $e$ positions in these $m$-tuples where all the components are simultaneously unbounded.  If this is not true then we can find some upper bound $B$ with the property that, if $(n_1,n_2,\ldots,n_m)\in Y_n$, then, for all subsets $E$ of size $e$ of the $m$ positions,  $n_i<B$ for some $i\in E$; but then it follows that $n_k<B$ for all but $e-1$ of the $m$ positions.  Hence the permutations of length $n$ in $X_{\pi}$ fall into $\binom{m}{e-1}$ subsets; in each subset there is an associated set of $e-1$ positions and outside of these positions the entries in the $m$-tuples are less than $B$.  Now it follows that each of the subsets has only $O(n^{e-2})$ elements, and this contradicts $\degree(X)=e-1$.

Therefore the sizes of maximal expansible sets determine $\degree(X)$ exactly.

The expansible sets can be characterised by an avoidance condition.  Suppose first that neither $p$ nor $q$ is zero.  Then, in order that no subpermutation isomorphic to $\beta$ appears when an expansible subset $E$ of an irreducible permutation $\theta$ is expanded, there must not be a subsequence $dbca$ of $\theta$ isomorphic to $4231$ with both $d,a\in E$.  Of course subpermutations isomorphic to $\alpha_r$ cannot arise by expansion at all.  In the case $p=0$ (or $q=0$) the condition is slightly different: we require that there should not be a subsequence $bca$ isomorphic to $231$ with $a\in E$ (or a subsequence $cab$ isomorphic to $312$ with $c\in E$).

Now, from these remarks and the upper bound on $\pi$ in the above lemma we have

\begin{proposition}
\begin{displaymath}\degree(\av(\alpha_r, \beta_{pq})) \leq
\left\{ \begin{array}{rl}
(r-1)^2(s-2) - r&\mbox{if $p > 0$ and $q > 0$,}\\
(r-1)^2(s-2)-1& \mbox{if $p = 0$ or $q = 0$}.
\end{array}\right.
\end{displaymath}
\end{proposition}

Unfortunately, we do not know how good this upper bound is.  A lower bound can be obtained as follows.

\begin{lemma}\label{lowerboundirreducible}
There exists an irreducible permutation of size $(r-1)(2s-5)$ in $\av(\alpha_{r},\beta_{pq})$; furthermore this permutation has an expansible set of size $(r-1)(s-2)$.
\end{lemma}

\begin{proof}
The  irreducible we construct consists of $r-1$ interlocking layers of decreasing points, each of size $2(p+q)-1$. For each layer, excluding the lower-leftmost, the bottom $p+q-1$ points are placed horizontally between the top $p+q$ points in the layer immediately below and to the left. Similarly, the top $p+q-1$ points are placed vertically in the gaps between the bottom $p+q$ points of the layer immediately below and to the left.

Thus every pair of consecutive points within a layer is separated by a point from the layer above and one from the layer below. This gives the permutation:

\begin{eqnarray*}
&&A_{1}+3x, A_{1}+3x-2,\ldots,  A_{1}+x, \\
&&A_{2}+4x, A_{1}+x-1, A_{2}+4x-2,\ldots, A_{2}+2x+2,A_{1}, \\
&&A_{2}+2x, A_{3}+4x , A_{2}+2x-2,A_{3}+4x-2,\ldots, A_{2},\\
&&A_{3}+2x, A_{4}+4x,\ldots,A_{3},\\
&&\ldots\\
&&A_{r-2}+2x, A_{r-1}+3x, A_{r-2}+2x-2, A_{r-1}+3x-1,\ldots,A_{r-2},\\
&&A_{r-1}+2x, A_{r-1}+2x-2,\ldots, A_{r-1}\\
\end{eqnarray*}
where $x=p+q-1$, and $A_{1} = 1$, $A_{2} = A_{1}+x+1$, $A_{i} = A_{i-1} + 2x+1$ for $3\leq i \leq r-1$ are the symbols at the bottom right of each layer (counting layers from bottom left to top right).

Figure \ref{lowerbound} shows the layers as a series of interlocking boxes, with each adjacent pair of layers shown in greater detail.

\begin{figure}[ht]
\begin{center}
\includegraphics[width=3in]{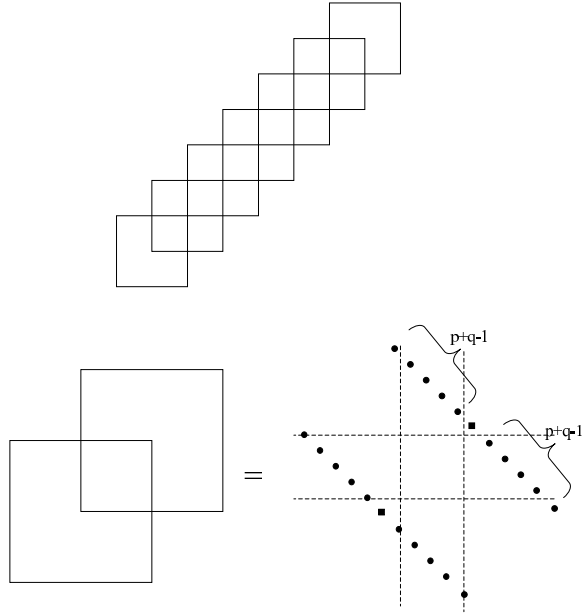}
\caption{An irreducible permutation of length $(r-1)(2s-5)$ in $\av(\alpha_{r},\beta_{pq})$}
\label{lowerbound}
\end{center}
\end{figure}

It is clear that the permutation above is irreducible, and that it does not contain any occurrence of $\alpha_{r}$ (as there are only $r-1$ layers).

If it were to contain $\beta_{pq}$, then the occurrence of $\beta_{pq}$ could only lie in 2 adjacent layers, as two non-adjacent layers do not intersect. Thus the top descent of length $p$ of such a $\beta_{pq}$ in the irreducible could be taken to lie at the top of the upper layer (or the top of the lower layer, respectively), and the bottom descent of length $q$ at the bottom of the lower layer (the bottom of the upper layer, respectively). But in either case this then leaves no points in between for the required increasing layer of size $2$ within $\beta$. Hence this is an irreducible permutation lying in $\av(\alpha_r, \beta_{pq})$, and it is of size $(r-1)(2(p+q)-1) = (r-1)(2s-5)$.

Now we identify an expansible subset of this permutation.
Suppose first that $p>0$ and $q>0$. We claim that the subset consisting only of the bottom $p+q$ points from each layer is expansible.

Note that such a subpermutation consists of exactly $r-1$ non-intersecting layers, each of size $p+q$. None of the bottom $p+q$ points of any layer can be used as the top point in a 4231-pattern, as there are no increases below and to the right of any one of these points. Thus the subpermutation is expansible, and is of size $(r-1)(p+q) = (r-1)(s-2)$.

In the case $p=0$ the same proof can be used, but we need to use the top $p+q$ points of each layer instead so that none can be used as a 1 in a $231$-pattern.
\end{proof}

\begin{corollary}$\degree(\av(\alpha_r, \beta_{pq}))\geq (r-1)(s-2)-1$.
\end{corollary}

\subsection{Case {\boldmath $s=3$}}

If $s=3$ then $\beta=231$ or $312$.  To within a symmetry that preserves $\alpha$ these are the same so we take $\beta=231$.  Now put
\[F_r(x)=\sum_{n=0}^{\infty}f_{rn}x^n\]
which is the generating function of the enumeration sequence $(f_{rn})_{n=0}^{\infty}$ for $\av(\alpha_r,\beta)$.  We now appeal to a result of Mansour and Vainshtein (Theorem 2.1 of \cite{MV}).  It tells us that
\[
F_r=1+x\sum_{i=1}^{r-1}(F_i-F_{i-1})F_{r+1-i}.
\]
Since both $\alpha_r$ and $\beta$ are irreducible, the class $\av(\alpha_r, \beta)$ can be obtained from its irreducible permutations by arbitrarily expanding individual elements into descending segments. Such classes were investigated in \cite{wreath} where it was shown that, if we define $G_r(y)=F_r(y/(1+y))$, then $G_r$ will be the generating function for the irreducible elements of $\av(\alpha_r, \beta)$.  After some manipulation we obtain
\[
G_r=1+y+y\sum_{i=2}^{r-1}(G_i-G_{i-1})G_{r+1-i}
\]
Although this equation is similar to the previous equation it differs in that the summation starts at $i=2$ and this means that $G_r$ does not appear on the right-hand side of the equation.  It follows that each $G_r$ is a polynomial in $y$.

\begin{lemma}
$G_r$ has degree $2r-3$ and leading coefficient the $(r-2)$th Catalan number $\cat(r-2)$ for all $r\geq 2$.
\end{lemma}
\begin{proof}
We have $G_0=0, G_1=1$ and all other $G_r$ are given by the recurrence.  Let $d_r$ be the degree of $G_r$.  Then we have $d_1=0$ and, from the recurrence,
\[
d_r=1+\max_{2\leq i\leq r-1}(d_i+d_{r+1-i})
\]
and $d_r=2r-3$ follows by induction.

Next, if $\lambda_r$ is the leading coefficient of $G_r$ we obtain
\[
\lambda_r=\sum_{i=2}^{r-1}\lambda_i\lambda_{r+1-i}
\]
and, again by induction, $\lambda_r=\cat(r-2)$.
\end{proof}

\begin{proposition}
The pattern class $\av(12\cdots r,231)$ is enumerated by a polynomial of degree $2r-4$ with leading coefficient $1 / (r-1)!(r-2)!$.
\end{proposition}
\begin{proof}
As noted above, this follows from the results of \cite{wreath}, specifically that
\[
F_r(x) = G_r \left( \frac{x}{1-x} \right).
\]
This equation simply captures symbolically the fact that each permutation belonging to $\av(\alpha_r, \beta)$ is obtained uniquely from the expansion of some  irreducible in the class.
\end{proof}

\subsection{Case {\boldmath $r=3$}}
In this subsection we shall take $\alpha=123$.

\begin{lemma}
Let $\gamma$ be an arbitrary irreducible permutation belonging to $\av(123,\beta_{pq})$. If $\delta$ is an expansible subset of
$\gamma$, then $|\delta| \leq 2(p+q)+2$ when $p>0$ and $q>0$, and $|\delta|\leq 2(p+q)+1$ when $p=0$ or $q=0$.
\end{lemma}
\begin{proof}
First assume $p>0$ and $q>0$, and note that an expansible subset of $\gamma$ will necessarily avoid 4231. Thus, $\delta$  avoids $123$ and $4231$.

Let $k=|\delta|$.  By Proposition 3.2 of \cite{restricted},  $\delta = \delta_{1}\delta_{2}\cdots\delta_{6}$, where, in the  permutation $\bar{\delta} = \bar{\delta_{1}}\bar{\delta_{2}}\cdots\bar{\delta_{6}}$ that is equivalent to $\delta$, each $\bar{\delta_i}$ is consecutive decreasing. Thus in $\bar{\delta}$ there are at least $k-6$ consecutive decreasing pairs $i+1,i$.

Each of these, as  pairs of $\gamma$, must be separated by some point of $\gamma$.  We shall show that two pairs $(f,e), (v,u)$ cannot be separated by the same term of $\gamma$. Assume without loss of generality that $(f,e)$ lies to the left of $(v,u)$. There are two cases:

\begin{enumerate}
\item $e<f<u<v$: here separating both pairs would create a 123 pattern.
\item $u<v<e<f$: here separating both pairs would create a $4231$-pattern, contradicting $\delta$ being expansible.
 \end{enumerate}

Hence $|\gamma| \geq k + (k-6) = 2k -6$, but (by Lemma \ref{maxdescent}) $|\gamma|\leq 4(p+q)-2$ so $k\leq 2(p+q)+2$.

In the case when $p=0$, we know that an expansible subset must now avoid $231$. So $\delta$ is easily seen to have $3$ segments rather than $6$ and so $2k-3\leq |\gamma|\leq 4(p+q)$, giving $k \leq 2(p+q)+1$.
\end{proof}

This lemma gave an upper bound on the size of an expansible set.  A matching lower bound can be obtained from the construction in Lemma \ref{lowerboundirreducible} specialised to $r=3$.  In the permutation given there we simply take the bottom $p+q+1$ points from each of the two decreasing sequences defined, and verify that this set of $2(p+q+1)$ points is expansible.  Together with the last lemma this proves the first statement of the next proposition; the last statement is a small variation on the first whose details we omit.

\begin{proposition}
If $p>0$ and $q>0$ then $\degree(123, \beta_{pq})=2s-3$.  If either $p=0$ or $q=0$, $\degree(123, \beta_{pq})=2s-4$.
\end{proposition}

\section{Discussion and open problems}
We have given a criterion for a pattern class to have polynomial
growth in terms of its basis and this enabled us to give complete
sets of two or three pattern restrictions to produce polynomial
growth.  However, the criterion produces neither the enumerating
polynomial nor its degree, and so we gave some properties of these
polynomials in the case of two restrictions.  In general, it remains
open to determine exactly what the polynomials are. This work might
be regarded as a first step in characterising possible enumeration
functions for a pattern class.

The authors wish to thank the referees for the improvements brought
about from their helpful comments and suggestions, and also Vince Vatter for his input.

\end{document}